\documentclass{article}
\usepackage{amsmath}[1996/11/01]
\usepackage{fullpage}
\usepackage{amssymb,amsthm,amsxtra}
\usepackage{epsfig}

\newtheorem{thm}{Theorem}

\newtheorem{lem}[thm]{Lemma}
\newtheorem{prop}[thm]{Proposition}
\theoremstyle{remark}
\newtheorem*{rem}{Remark}

\newtheorem*{notrem}{Notational Remark}

\theoremstyle{definition}

\numberwithin{equation}{section}
\numberwithin{thm}{section}

\DeclareMathOperator{\supp}{supp}

\begin{document}

\title{\bf Eigenvalues of Large Sample Covariance Matrices of Spiked
Population Models}
\author{{\bf Jinho Baik}\footnote{Department of Mathematics, University
of Michigan, Ann
Arbor, MI, 48109, USA, baik@umich.edu}\ \ and {\bf Jack W.
Silverstein}\footnote{Department of Mathematics, North Carolina
State University, Raleigh, NC, 27695, USA,
jack@math.ncsu.edu
\break \hskip.3in AMS 1991 {\sl subject classifications}. Primary 15A52,
60F15; secondary 62H99.\hfil\break
}}

\date{July 27, 2004}
\maketitle

\begin{abstract}
We consider a spiked population model, proposed by Johnstone, whose
population eigenvalues are all unit except for a few fixed
eigenvalues. The question is to determine how the sample
eigenvalues depend on the non-unit population ones
when both sample size and population size become large. This paper
completely determines the almost sure limits for a general class
of samples.
\end{abstract}

\section{Introduction}

The sample covariance matrix is fundamental to multivariate
statistics. When the population size is not large and for a
sufficient number of samples, the sample covariance matrix is a
good approximate of the population covariance matrix . However
when the population size is large and comparable with the sample
size, as is in many contemporary data, it is known that the sample
covariance matrix is no longer a good approximation to the
covariance matrix. The Marchenko-Pastur theorem \cite{MP} states
that with $n=$ the sample size, $p=$ the population size, as
$n=n(p) \to \infty$ such that $\frac{p}{n} \to c$, the eigenvalues
$s_j^{(p)}$, $j=1, \dots, p$, of the sample covariance matrix of
normalized i.i.d. Gaussian samples satisfy for any real $x$
\begin{equation}\label{eq:DOS}
  \frac1{p} \#\{s_j^{(p)}: s_j^{(p)} <x\} \to F(x)
\end{equation}
almost surely where
\begin{equation}\label{eq:PM}
   F'(x) = \frac{1}{2\pi xc}\sqrt{(b-x)(x-a)}, \qquad a<x<b,
\end{equation}
and $a=(1-\sqrt{c})^2$ and $b=(1+\sqrt{c})^2$ when $0< c \le 1$.
When $c>1$, there is an additional Dirac measure at $x=0$ of mass
$1-\frac1{c}$. Moreover, there are no stray eigenvalues in the
sense that the top and bottom eigenvalues converge to the edges of
the support of $F$ \cite{Geman}:
\begin{equation}\label{eq:s1}
  s_1^{(p)} \to (1+\sqrt{c})^2
\end{equation}
almost surely and \cite{Silverstein85}
\begin{equation}\label{eq:smin}
  s_{\min\{p, n\}}^{(p)} \to (1-\sqrt{c})^2
\end{equation}
almost surely ($s^{(p)}_{n+1}=\dots s^{(p)}_p=0$ when $n<p$). One
can extract from this some information of the population
covariance matrix even though the sample covariance matrix is not
a good approximate. For example, if there are non-zero eigenvalues
of the sample covariance matrix well separated from the rest of the
eigenvalues, one finds, assuming the Gaussian entries, that the
samples are not i.i.d..

There are indeed many cases in which a few eigenvalues of the
sample covariance matrix are separated from the rest of the
eigenvalues, the latter being packed together as in the support of the
Marchenko-Pastur function \eqref{eq:PM}. The examples include
speech recognition \cite{Buja, Johnstone} , mathematical finance
\cite{PlerousGRAGS}, \cite{LalouxCPB}, \cite{MalevergneS},
wireless communication \cite{Telatar}, physics of mixture
\cite{SearC}, and data analysis and statistical learning
\cite{HoyleR}.

As a possible explanation for such features, Johnstone
\cite{Johnstone} proposed the `spiked population model' where all
but finitely many eigenvalues of the population covariance matrix
are the same, say equal to $1$. The question is how the
eigenvalues of the sample covariance matrix would depend on the
non-unit population eigenvalues as $p, n \to\infty$. It is
known \cite{MP, SB95} that the Marchenko-Pastur result
\eqref{eq:DOS} still holds for the spiked model. But \eqref{eq:s1}
and \eqref{eq:smin} are not guaranteed and some of the eigenvalues
are not necessarily in the support of \eqref{eq:PM}.

For example, consider the case when the population covariance
matrix has one non-unit eigenvalue, denoted by $\sigma_1$. If
$\sigma_1$ is close to $1$, one would expect that as the dimension
$p$ becomes large, the population covariance matrix would be close
to a large identity matrix, and hence $\sigma_1$ would have little
effect on the eigenvalues of the sample covariance matrix. On the
other hand, if $\sigma_1$ is much bigger than $1$, then even if
$p$ becomes large, $\sigma_1$ might still pull up the eigenvalues
of the sample covariance matrix. How big should $\sigma_1$ be in
order to have any effect, how many eigenvalues of the sample
covariance matrix would be pulled up and exactly where would the
pulled-up eigenvalues be ? We will see in the results below that the
answers are $\sigma_1>1+\sqrt{c}$ (where $\frac{p}{n}\to c$), one
eigenvalue at most, and $\sigma_1+\frac{c\sigma_1}{\sigma_1-1}$,
respectively.

For \emph{complex Gauassian} samples, the papers \cite{Sandrine,
BBP} study the \emph{largest} eigenvalue of the sample
covariance matrix. The authors determine the transition behavior
and the limiting distributions are also obtained. The purpose of
this paper is a complete study of the spiked model for \emph{both
real and complex samples which are not necessarily Gaussian}. We
obtain almost sure limit results. A general study of `non-null'
covariance matrices was done in \cite{BaiS98, BaiS99}. We will
show in this paper how to extract the desired results from the
work of \cite{BaiS99}. While this paper was being prepared, the
authors learned that Debashis Paul \cite{DPaul:unp} was also
studying the spiked model independently at the same time, which
has some overlap with this work. Paul considers the real Gaussian
samples for $c<1$, and obtains the almost sure limits as in
\eqref{eq:thm11} and \eqref{eq:thm12} below for large sample
eigenvalues. Moreover, when all non-unit population eigenvalues
are simple, the limiting distribution is found to be Gaussian (see
Subsection \ref{sec:discussion} below for more detail). On the
other hand, our paper (i) is concerned with more general samples,
not necessarily Gaussian, (ii) includes all choices of $c$ and
(iii) studies both large and small sample eigenvalues. We remark
that a complete study of limiting distributions is still an open
question.


\subsection{Model}

Let $T_p$ be a fixed $p\times p$ non-negative definite Hermitian matrix.
Let $Z_{ij}$, $i,j=1,2,\dots$, be independent and identically distributed
complex valued random variables satisfying
\begin{equation}
   \mathbb{E}(Z_{11})=0, \qquad
\mathbb{E}(|Z_{11}|^2)=1,
\quad \text{and} \quad \mathbb{E}(|Z_{11}|^4)<\infty,
\end{equation}
and set $Z_p=(Z_{ij})$, $1\le i\le p$, $1\le j\le n$. We take the
sampled vectors to be the columns of $T_p^{1/2}Z_p$, where $T_p^{1/2}Z_p$
is an Hermitian square root of $T_p$. Hence $T_p$ is the population
covariance matrix.  Of course, not all random vectors are realized
as such, but this model is still very general. When $Z_{ij}$ are
i.i.d (real or complex) Gaussian, the model becomes the Gaussian samples with
population covariance matrix $T_p$. Outside the Gaussian case we see that
these vectors cover a broad range of random vectors, completely real or
complex, with arbitrary population covariance matrix.


Let
\begin{equation}
  B_p:= \frac1{n} T_p^{1/2} Z_p Z_p' T_p^{1/2}
\end{equation}
be the sample covariance matrix, where $Z_p'$ denotes conjugate transpose.
Denote the eigenvalues of $B_p$ by $s_1^{(p)}, \dots, s_p^{(p)}$: for some
unitary matrix $U_B$,
\begin{equation}
  U_BB_pU_B^{-1} =
  \begin{pmatrix}
  s_1^{(p)} \\ & s_2^{(p)} \\ & & \ddots \\ & & & s^{(p)}_p
  \end{pmatrix}
  = diag( s_1^{(p)}, s_2^{(p)}, \dots, s^{(p)}_p).
\end{equation}
For definiteness, we order the eigenvalues as $s_1^{(p)}\ge
s_2^{(p)}\ge \dots\ge s^{(p)}_p\ge 0$.

Let $\alpha_1>\dots> \alpha_M>0$ be fixed real numbers for some
fixed $M\ge 0$, which is independent of $p$ and $n$. Let $k_1,
\dots, k_{M}$ be fixed non-negative integers and set $r=k_1+\dots
+ k_M$, which are also independent of $p$ and $n$. We assume that
all the eigenvalues of $T_p$ are $1$ except for, say, the first
$r$ eigenvalues. This is the `spiked population model' proposed in
\cite{Johnstone}. Let the first $r$ eigenvalues be equal to
$\alpha_1, \dots, \alpha_M$ with multiplicity $k_1, \dots, k_M$,
respectively: for some unitary matrix $U_T$,
\begin{equation}\label{eq:Tpdef}
\begin{split}
  U_TT_pU_T^{-1}
  &= diag(
  \underbrace{\alpha_1, \dots, \alpha_1}_{k_1}, \underbrace{\alpha_2, \dots,
\alpha_2}_{k_2},
  \dots, \underbrace{\alpha_M, \dots, \alpha_M}_{k_M},
  \underbrace{1, \dots, 1}_{p-r} ) .
\end{split}
\end{equation}
We set $k_0=0$.

\subsection{Results}\label{sec:results}

\begin{thm}[case $c<1$]\label{thm1}
Assume that $n=n(p)$ and $p\to\infty$ such that
\begin{equation}
  \frac{p}{n} \to c
\end{equation}
for a constant $0<c<1$. Let $M_0$ be the number of $j$'s such that
$\alpha_j>1+\sqrt{c}$, and let $M-M_1$ be the number of $j$'s such
that $\alpha_j<1-\sqrt{c}$.
Then the following holds.
\begin{itemize}
\item For each $1\le j\le M_0$,
\begin{equation}\label{eq:thm11}
  s^{(p)}_{k_1+\dots + k_{j-1}+i} \to \alpha_{j} + \frac{c
  \alpha_{j}}{\alpha_j-1},
  \qquad 1\le i\le k_j.
\end{equation}
almost surely.
\item
\begin{equation}\label{eq:thm12}
  s^{(p)}_{k_1+\dots+k_{M_0}+1} \to (1+\sqrt{c})^2
\end{equation}
almost surely.
\item
\begin{equation}
  s^{(p)}_{p-r+k_1+\dots + k_{M_1}} \to (1-\sqrt{c})^2
\end{equation}
almost surely (recall $r=k_1+\dots + k_M$).
\item
For each $M_1+1\le j\le M$,
\begin{equation}\label{eq:thm13}
  s^{(p)}_{p-r+k_1+\dots + k_{j-1} +i}\to  \alpha_{j} + \frac{c
  \alpha_{j}}{\alpha_j-1},
  \qquad 1\le  i\le k_{j}
\end{equation}
almost surely.
\end{itemize}
\end{thm}

Therefore, when $c<1$, in order for a population eigenvalue to
contribute a non-trivial effect to the eigenvalues of the sample
covariance matrix, it should sufficiently big (larger than
$1+\sqrt{c}$) or sufficiently small (less than $1-\sqrt{c}$). As
an example, when $r=1$, by denoting the only non-unit eigenvalue
by $\sigma_1$, the largest sample eigenvalue $s^{(p)}_1$ satisfies
\begin{equation}
  s^{(p)}_1 \to \begin{cases}
  (1+\sqrt{c})^2, \qquad & \sigma_1\le 1+\sqrt{c} \\
  \sigma_1 + \frac{c\sigma_1}{\sigma_1-1}, \qquad & \sigma_1 >
  1+\sqrt{c}
  \end{cases}
\end{equation}
almost surely. When $r=2$, by denoting the two non-unit
eigenvalues by $\sigma_1, \sigma_2$, the largest sample eigenvalue
$s^{(p)}_1$ satisfies
\begin{equation}
  s^{(p)}_1 \to \begin{cases}
  (1+\sqrt{c})^2, \qquad & \max\{\sigma_1, \sigma_2\} \le 1+\sqrt{c} \\
  \max\{\sigma_1, \sigma_2\} + \frac{c\max\{\sigma_1, \sigma_2\}}{\max\{\sigma_1,
  \sigma_2\}-1},
  \qquad & \max\{\sigma_1, \sigma_2\} > 1+\sqrt{c}
  \end{cases}
\end{equation}
almost surely.

The results \eqref{eq:thm11} and \eqref{eq:thm12} are also
independently obtained in \cite{DPaul:unp} under the assumption
that the samples are Gaussian.

\begin{thm}[case $c>1$]\label{thm2}
Assume that $n=n(p)$ and $p\to\infty$ such that
\begin{equation}
  \frac{p}{n} \to c
\end{equation}
for a constant $c>1$. Let $M_0$ be the number of $j$'s such that
$\alpha_j>1+\sqrt{c}$. Then the following holds.
\begin{itemize}
\item
For each $1\le j\le M_0$,
\begin{equation}\label{eq:thm21}
  s^{(p)}_{k_1+\dots + k_{j-1}+i} \to \alpha_{j} + \frac{c
  \alpha_{j}}{\alpha_j-1},
  \qquad 1\le i\le k_j.
\end{equation}
almost surely.
\item
\begin{equation}\label{eq:thm22}
  s^{(p)}_{k_1+\dots +k_{M_0}+1} \to (1+\sqrt{c})^2
\end{equation}
almost surely.
\item
\begin{equation}\label{eq:thm23}
   s^{(p)}_n \to (1-\sqrt{c})^2
\end{equation}
almost surely.
\item For all $p$,
\begin{equation}\label{eq:thm24}
   s^{(p)}_{n+1}=\dots = s^{(p)}_p=0.
\end{equation}
\end{itemize}
\end{thm}

Thus, unlike the case of $c<1$, small eigenvalues of $T_p$
do not affect the eigenvalues of $B_p$ when $c>1$.

\begin{thm}[case $c=1$]\label{thm3}
Assume that $n=n(p)$ and $p\to\infty$ such that
\begin{equation}
  \frac{p}{n} \to 1.
\end{equation}
Let $M_0$ be the number of $j$'s such that $\alpha_j>2$. Then the
following holds.
\begin{itemize}
\item
For each $1\le j\le M_0$,
\begin{equation}\label{eq:thm31}
  s^{(p)}_{k_1+\dots + k_{j-1}+i} \to \alpha_{j}
  + \frac{\alpha_{j}}{\alpha_j-1},
  \qquad 1\le i\le k_j.
\end{equation}
almost surely.
\item
\begin{equation}\label{eq:thm32}
  s^{(p)}_{k_1+\dots +k_{M_0}+1} \to 4
\end{equation}
almost surely.
\item
\begin{equation}\label{eq:thm33}
   s^{(p)}_{\min\{n,p\}} \to 0
\end{equation}
almost surely.
\end{itemize}
\end{thm}

\subsection{Discussion}\label{sec:discussion}

As mentioned earlier, the limiting density of the eigenvalues of
spiked population models is given by the Marchenko-Pastur theorem
\eqref{eq:PM} as in the identity population matrix case, and for
the top eigenvalue $s^{(p)}_1$ in the \emph{complex Gaussian}
case, the results Theorem 1.1 and Theorem 1.3 were first obtained
in \cite{Sandrine, BBP}. The paper \cite{BBP} (see section 6)
contains an interesting heuristic argument for the critical value
$1+\sqrt{c}$ and the value
\begin{equation}\label{eq:alphamysterious}
   \alpha_j + \frac{c\alpha_j}{\alpha_j-1}
\end{equation}
for $1\le j\le M_0$ above: they come from a competition between a
1-dimensional last passage time and a 2-dimensional last passage
time. It would be interesting to have such a heuristic reasoning
for the general case.

When $T_p$ is the identity matrix (the `null case'), under the
Gaussian assumption, the limiting distribution for the largest
eigenvalue is obtained for the complex case in \cite{Forrester,
kurtj:disc} and for the real case in \cite{Johnstone}.
\cite{Soshnikovcovariance} shows the Gaussian assumption is not
necessary when $c=1$. The limiting distributions are the
Tracy-Widom distributions \cite{TracyWidom, TracyWidomGOE} in the
random matrix theory in mathematical physics. For the spiked model
with complex Gaussian samples when $c\le 1$, the limiting
distributions of the largest eigenvalue are obtained in
\cite{Sandrine, BBP}. The paper \cite{BBP} determines the limiting
distribution of $s^{(p)}_1$ for complete choices of the largest
population eigenvalue $\alpha_1$ and its multiplicity $k_1$: the
distribution is (i) the Tracy-Widom distribution when
$\alpha_1<1+\sqrt{c}$, (ii) certain generalizations of the
Tracy-Widom distribution (see also \cite{Baik}) when $\alpha_1 =
1+\sqrt{c}$, and (iii) the Gaussian distribution ($k_1=1$) and its
generalization ($k_1\ge 2$, the Gaussian unitary ensemble) when
$\alpha_1>1+\sqrt{c}$. For real Gaussian samples \cite{DPaul:unp}
showed that when $c<1$, $M_0\ge 1$ and $k_1=\dots=k_{M_0}=1$, the
limiting distribution of $s^{(p)}_j$, $1\le j\le M_0$, is
Gaussian. It is an interesting open question to determine the
limiting distribution for the general case of real samples. See
section 1.3 of \cite{BBP} for a conjecture for the scaling.


\begin{figure}[ht]
\centerline{\epsfxsize=9cm\epsfbox{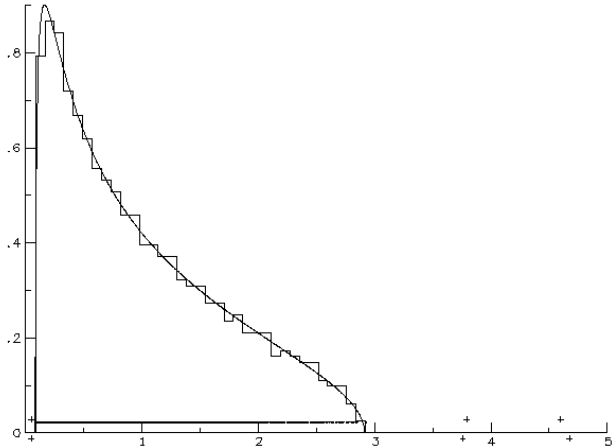}} \caption{Gaussian
samples when $p=1000, n=2000$}\label{fig:gauss2}
\end{figure}
\begin{figure}[ht]
\centerline{\epsfxsize=9cm\epsfbox{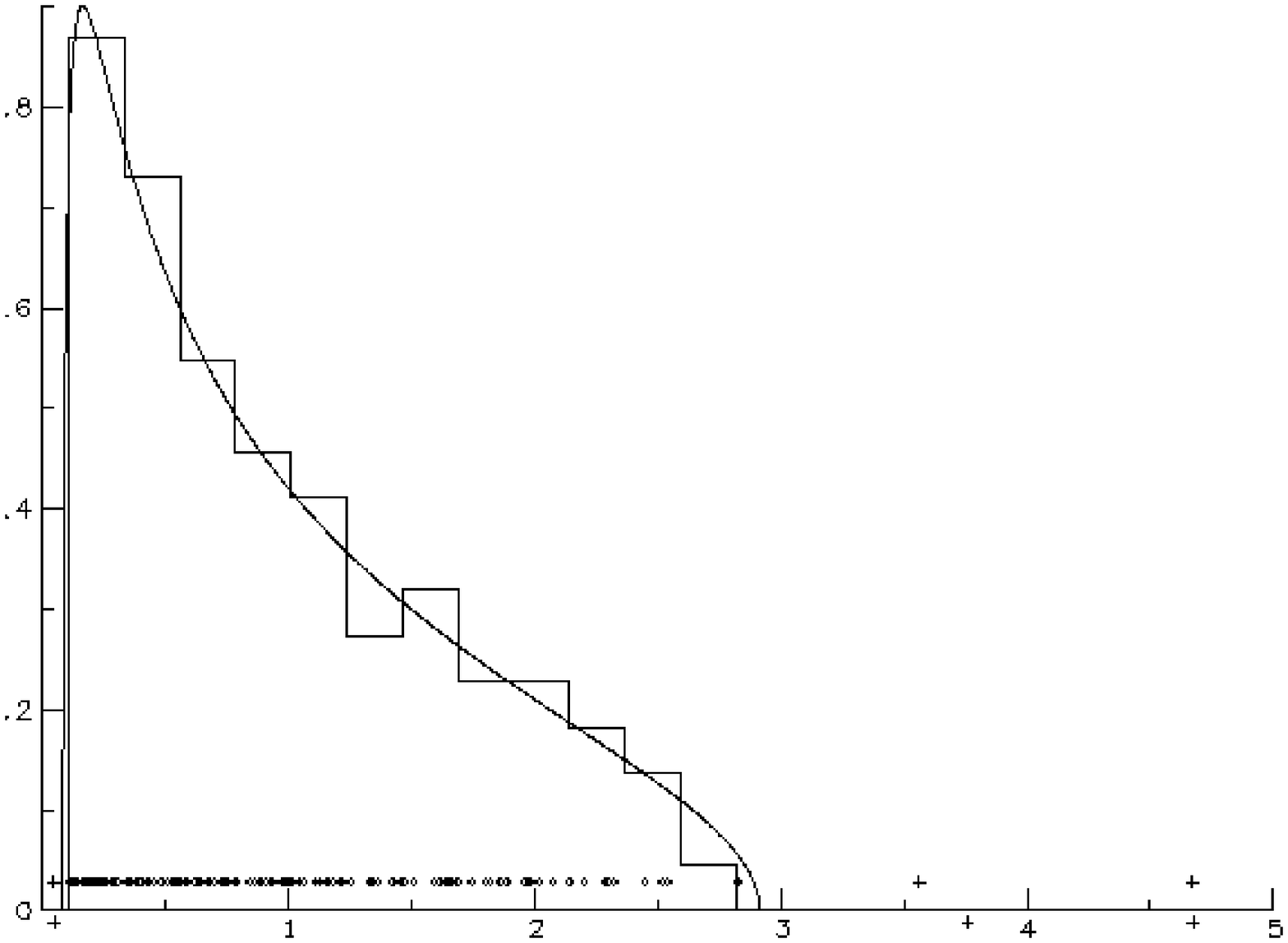}} \caption{Gaussian
samples when $p=100, n=200$}\label{fig:gauss3}
\end{figure}
\begin{figure}[ht]
\centerline{\epsfxsize=9cm\epsfbox{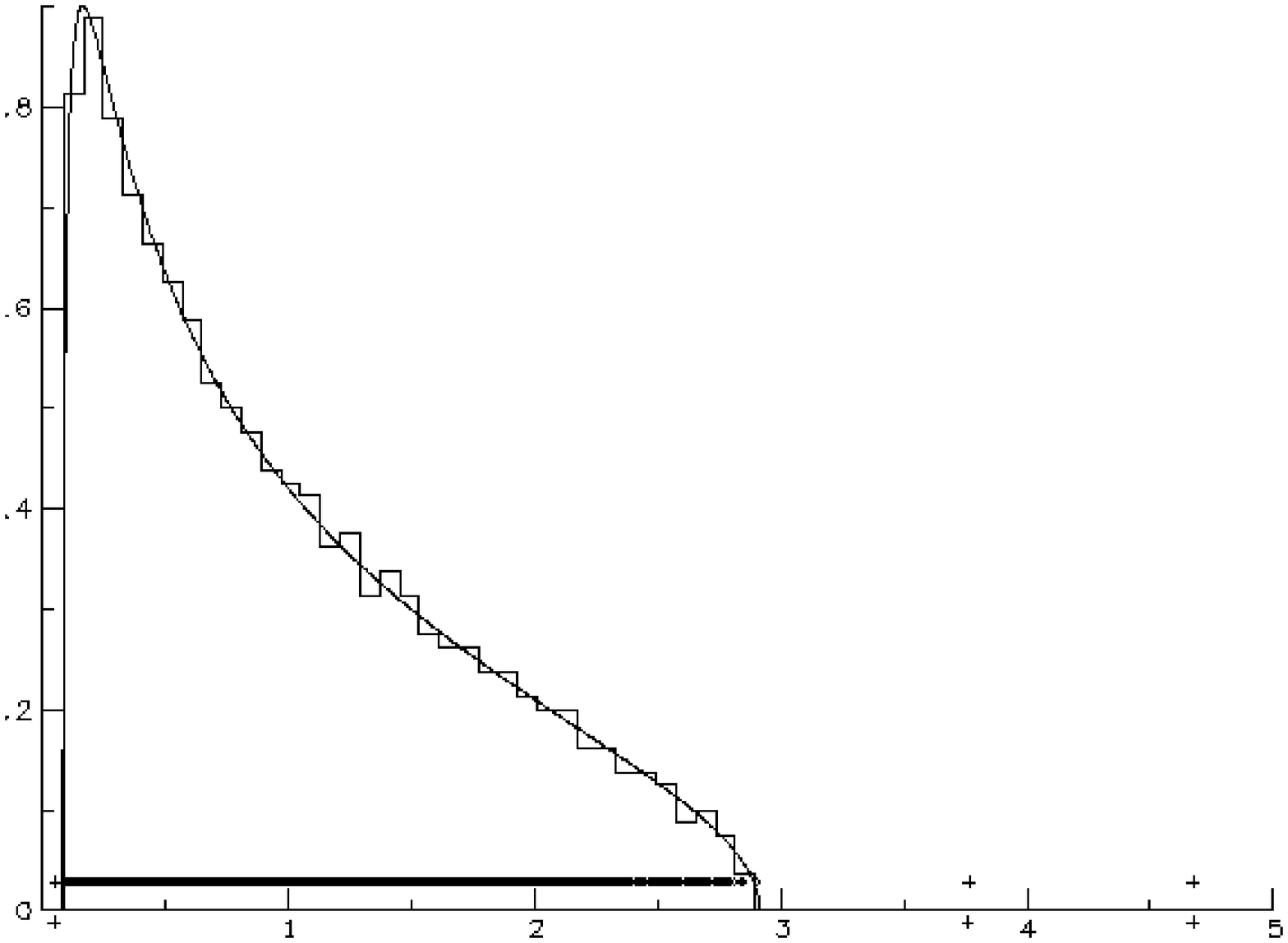}} \caption{Bernoulli
samples taking values $-1$ or $1$ when $p=1000,
n=2000$}\label{fig:uni1}
\end{figure}
\begin{figure}[ht]
\centerline{\epsfxsize=9cm\epsfbox{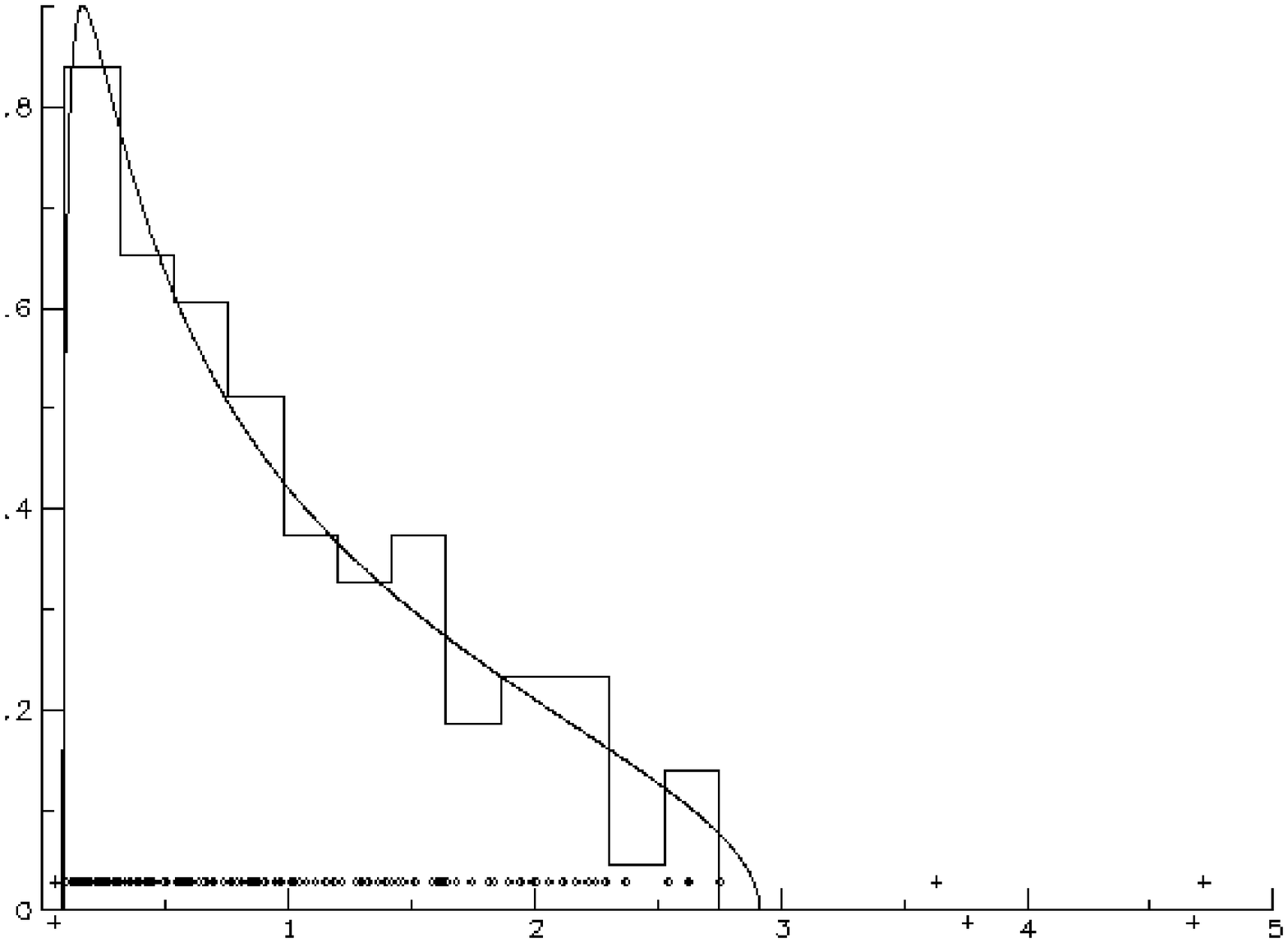}} \caption{Bernoulli
samples taking values $-1$ or $1$ when $p=100,
n=200$}\label{fig:uni100}
\end{figure}
We include several plots for the case when $c=0.5$ and there are
three non-unit population eigenvalues given by $0.1$, $3$ and $4$
(of multiplicity $1$ each). In this case, the critical values of
the eigenvalues are $1+\sqrt{c}\simeq 1.70711$ and $1-\sqrt{c}
\simeq 0.29289$. Hence theoretically we expect that three sample
covariance eigenvalues of values
$\alpha_j+\frac{c\alpha_j}{\alpha_j-1} \simeq 0.04444, 3.75$ and
$4.66667$ are away from the interval $[(1-\sqrt{c})^2,
(1+\sqrt{c})^2] \simeq [0.08578, 2.91422]$. The histogram and the
scatterplot of Figure \ref{fig:gauss2} is from Gaussian samples
when $p=1000, n=2000$. The smooth curve is the theoretical
limiting density and the theoretical locations of the three
separated eigenvalues are plotted with  $+$ signs below the
horizontal axis.  The smallest and largest two sample eigenvalues
are plotted with $+$ signs about the horizontal axis.  Figure
\ref{fig:gauss3} is from  Gaussian samples when $p=100, n=200$
while Figures \ref{fig:uni1} and \ref{fig:uni100} from samples of
Bernoulli variables taking values $-1$ or $1$ when $p=1000,
n=2000$, and $p=100, n=200$ respectively. The observed values of
the four separated eigenvalues in each case are as follows:

\bigskip
\begin{center}
\begin{tabular}{|c|c|c|c|}
  \hline  & smallest eigenvalue  & 2nd largest eigenvalue & largest eigenvalue \\ \hline
  theoretical & 0.04444 & 3.75 & 4.66667 \\ \hline
  Gaussian $p=1000$ & 0.04369 & 3.78400 & 4.59127 \\ \hline
  Gaussian $p=100$ &0.03979  & 3.55388 & 4.66192\\ \hline
  Bernoulli $p=1000$ & 0.04555 & 3.75706 & 4.66594 \\ \hline
  Bernoulli $p=100$ &0.05015 & 3.62337 & 4.70786 \\ \hline
\end{tabular}
\end{center}

\bigskip

\begin{figure}[ht]
\centerline{\epsfxsize=9cm\epsfbox{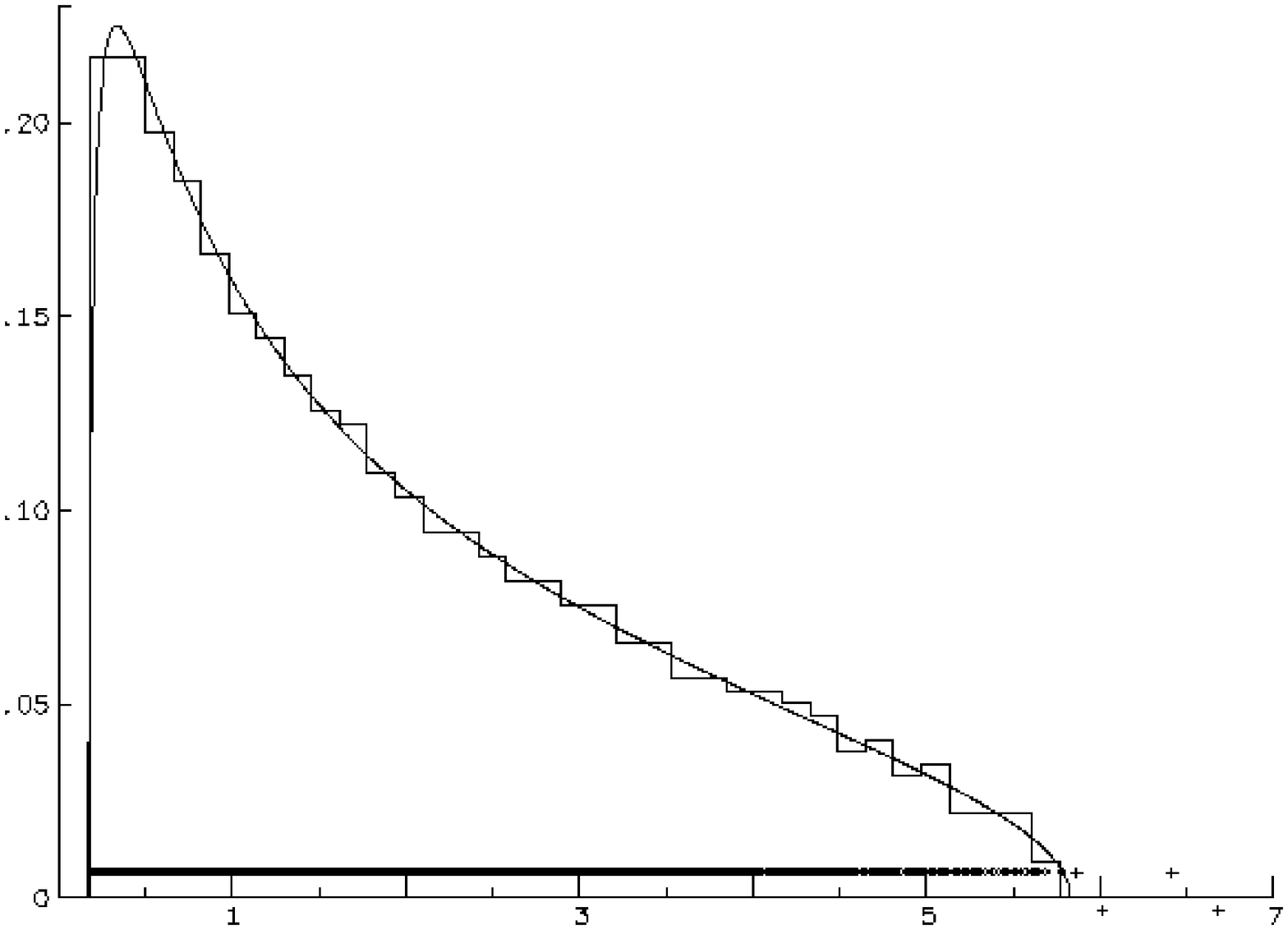}} \caption{Gaussian
samples when $p=2000, n=1000$}\label{fig:gauss4}
\end{figure}
\begin{figure}[ht]
\centerline{\epsfxsize=9cm\epsfbox{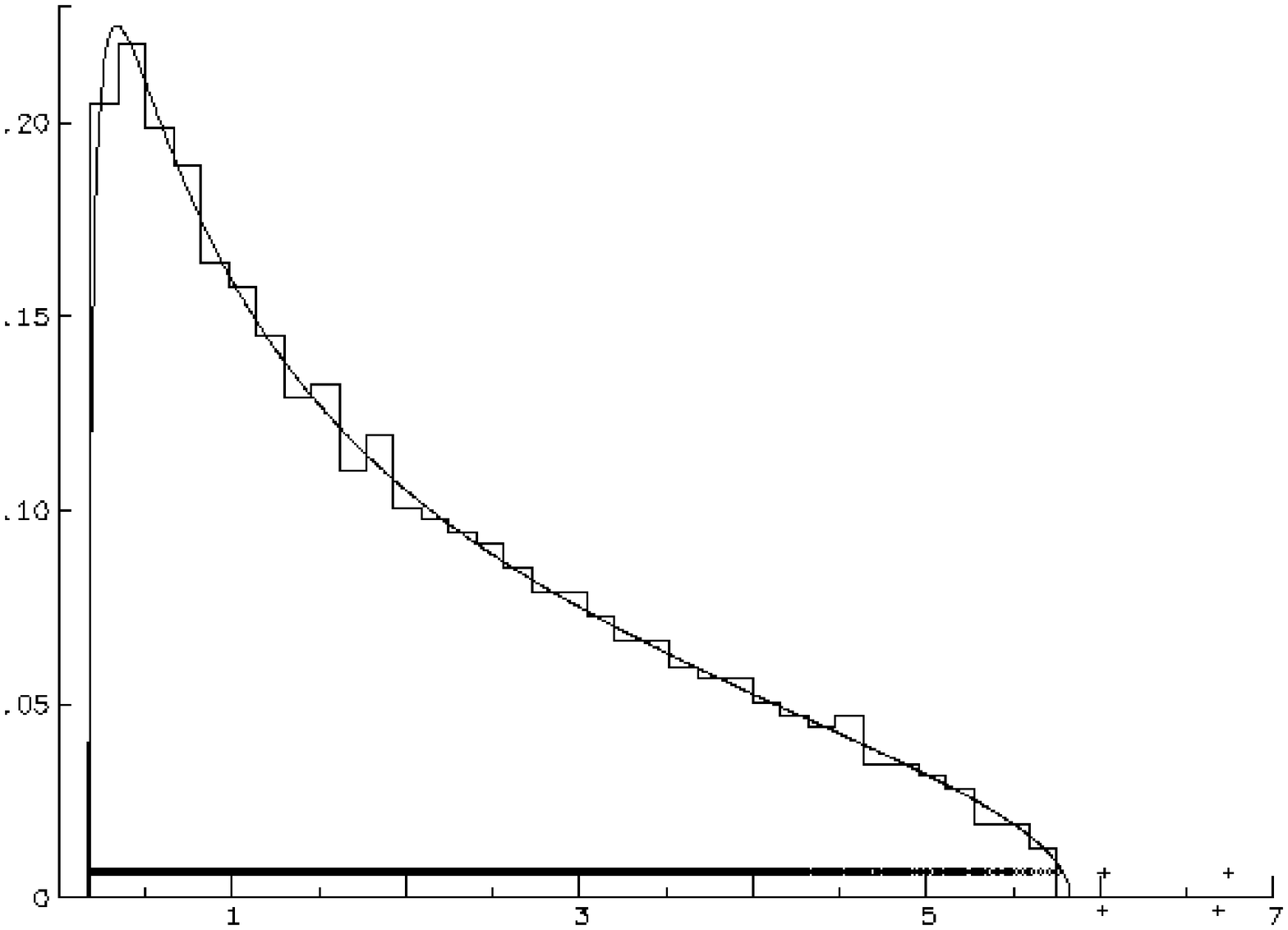}} \caption{Bernoulli
samples taking values $-1$ or $1$ when $p=2000,
n=1000$}\label{fig:uni2}
\end{figure}
Figure \ref{fig:gauss4} and Figure \ref{fig:uni2} are the cases
when $c=2$, $p=2000, n=1000$ with Gaussian and Bernoulli samples,
respectively. Again three non-unit population eigenvalues are
chosen $0.1$, $3$ and $4$. The critical value of the eigenvalues
is $1+\sqrt{c}\simeq 2.41421$ and the theory predicts that the two
largest sample eigenvalues given by
$\alpha_j+\frac{c\alpha_j}{\alpha_j-1} \simeq 6$ and $6.66667$ are
separated from the interval $[(1-\sqrt{c})^2, (1+\sqrt{c})^2]
\simeq [0.17157, 3.41209]$. Only non-zero eigenvalues are plotted
in Figure \ref{fig:gauss4} and Figure \ref{fig:uni2}. The observed
values of the separated eigenvalues in each case are as follows:

\bigskip
\begin{center}
\begin{tabular}{|c|c|c|}
  \hline  & 2nd largest eigenvalue & largest eigenvalue \\ \hline
  theoretical & 6 & 6.66667 \\ \hline
  Gaussian $p=2000$ & 5.8523 &  6.4013 \\ \hline
  Bernoulli $p=2000$ & 6.01065 & 6.725 \\ \hline
\end{tabular}
\end{center}
\bigskip

\bigskip

The paper is organized as follows. In section \ref{sec:BaiS}, we
summarize the work of Z. D. Bai and J. W. Silverstein on which we
heavily rely to prove our results. It turns out that the
determination of the support of a Stieltjes transform plays
the crucial role. This is obtained in section \ref{sec:support}.
The proofs of the main theorems are given in section
\ref{sec:proof}.

\noindent {\bf Acknowledgments.} Special thanks are due to Min Kang
for kindly inviting J.B. to give a talk at North Carolina State
University where the authors happened to have a chance to discuss
about the problem, which eventually lead to this work. We would
also like to thank Iain Johnstone for telling us the work of
Debashis Paul \cite{DPaul:unp} which was being done independently
and at the same time. The work of J.B. was supported in part by
NSF Grant \#DMS-0350729.

\section{Results of Z. D. Bai and J. W. Silverstein}\label{sec:BaiS}

Our analysis replies heavily on the work
\cite{BaiS98, BaiS99} of Bai and Silverstein.
In this section, we summarize the necessary results from \cite{BaiS98, BaiS99}.

\begin{notrem}
We denote by $p$ the population size and by $n$ the sample size.
The notations $n$ and $N$ are used in \cite{BaiS99} for $p$ and
$n$, respectively.
\end{notrem}

For a distribution function $G(\lambda)$, its
Stieltjes transform $m_G(z)$ is defined by
\begin{equation}
  m_G(z) = \int_{-\infty}^\infty \frac1{\lambda-z} dG(\lambda),
\qquad z\in \mathbb{C}^+:= \{ z\in\mathbb{C} : Im(z)>0\}.
\end{equation}
Also recall the inversion formula
\begin{equation}
   G([a,b]) = \frac1{\pi} \lim_{\eta\downarrow 0} \int_a^b
Im \bigl( m_G(\xi+i\eta) \bigr) d\xi
\end{equation}
for continuity points $a, b$ of $G$.


Assume the following:
\begin{itemize}
\item[(a)] $Z_{ij}$ are i.i.d. random variables in $\mathbb{C}$
with $\mathbb{E}(Z_{11})=0$, $\mathbb{E}|Z_{11}|^2=1$
and $\mathbb{E}|Z_{11}|^4<\infty$.
\item[(b)] $n=n(p)$ with $c_p:= \frac{p}{n} \to c>0$ as $p\to\infty$.
\item[(c)] For each $p$, $U_TT_pU_T^{-1}=diag(\sigma_1^{(p)}, \dots, \sigma_p^{(p)})$
for some unitary matrix $U_T$ such that
$H_p \to H_\infty$ in distribution for some distribution function $H_\infty$ where
$H_p$ is the empirical distribution function of the eigenvalues of $T_p$ defined by
\begin{equation}
  dH_p(\lambda) = \frac1{p} \sum_{j=1}^p \delta_{\sigma_j^{(p)}} (\lambda).
\end{equation}
\item[(d)] $\max\{ \sigma_1^{(p)}, \dots, \sigma_p^{(p)}\}$ is bounded in $p$.
\item[(e)] Set $Z_p=(Z_{ij})$, $1\le i\le p$, $1\le j\le n$
and $B_p= \frac1{n} T_p^{1/2}Z_pZ_p^* T_p^{1/2}$.
\item[(f)]
Set
\begin{equation}
  z_p(m) = -\frac1{m} + c_p \int \frac{t}{1+tm} dH_p(t).
\end{equation}
>From \cite{SB95} and \cite{Silverstein95},
it is known that there is a unique inverse function
$m_p(z)$ such that $m_p(z)\in \mathbb{C}^+$ for $z\in\mathbb{C}^+$
 It is also known \cite{SB95, Silverstein95}
that $m_p(z)$ is the Stieltjes transform
of a distribution, which will be denoted by $F_p$:
\begin{equation}\label{eq:defFpviamp}
  m_p(z) = \int_{-\infty}^\infty  \frac1{\lambda-z} dF_p(\lambda),
\qquad z\in \mathbb{C}^+.
\end{equation}
Suppose that the interval $[a,b]$ with $a>0$ lies in an open
interval outside the support of $F_p$ for all large $p$.
\end{itemize}

Now we will state the main result of \cite{BaiS99} which we need for our analysis.
It is easy to check that $F_p$ converges to some distribution
function $F_\infty$.
Then \cite{Silverstein95} $F_\infty(\lambda)$ is the almost sure
limit of the empirical spectral
distribution of $\underline{B}_p:=\frac1{n} Z_p^*T_pZ_p$ and
\begin{equation}\label{eq:Flam}
  F(\lambda):= \frac1{c} \bigl( F_\infty(\lambda) - (1-c) 1_{[0,\infty)} \bigr)
\end{equation}
is the almost sure limit of the empirical spectral distribution of
$B_p= \frac1{n} T_p^{1/2}Z_pZ_p^* T_p^{1/2}$.

\begin{rem}
The function $F_p$ is not the empirical distribution of $B_p$.
The distribution function $F_p$ is defined only through \eqref{eq:defFpviamp}.
\end{rem}

Moreover \cite{BaiS98}, the Stieltjes transform of $F_\infty$,
\begin{equation}\label{eq:minf}
  m_\infty(z) = \int_{-\infty}^\infty  \frac1{\lambda-z} dF_\infty(\lambda),
\qquad z\in \mathbb{C}_+:= \{ z\in\mathbb{C} : Im(z)>0\},
\end{equation}
is invertible, with the inverse given by
\begin{equation}\label{eq:zinf0}
  z_\infty(m) = -\frac1{m} + c \int \frac{t}{1+tm} dH_\infty(t).
\end{equation}
On the other hand, given $H_\infty$, $F_\infty$ is determined from
\eqref{eq:zinf0} and \eqref{eq:minf}.
Note that $m_\infty(z)$ is
well-defined not only on $\mathbb{C}_+$ but also up to the real
line outside $\supp(F_\infty)$ and its inverse exists on
$m_\infty(\mathbb{C}_+\cup \supp(F_\infty)^c)$.

\begin{rem}
If $[a,b]$ satisfies condition (f) above, it is easy to check that
$[a, b]\subset \supp(F_\infty)^c$.
\end{rem}

Given an interval $[a,b]$ satisfying condition (f) above
and $m_\infty(b)<0$, it is shown in \cite{BaiS99} that
there is an integer $i_p\ge 0$ satisfying the conditions
\begin{equation}\label{eq:mainpropeqcon}
  \sigma_{i_p}^{(p)} > -\frac1{m_\infty(b)},
\qquad \sigma_{i_p+1}^{(p)}<-\frac1{m_\infty(a)}
\end{equation}
for large $p$. (Here $\sigma^{(p)}_0:= \infty$.)

\begin{prop}[Theorem 1.2 \cite{BaiS99}]\label{prop:mainprop}
Assume (a)-(f) above.
\begin{itemize}
\item[(i)] If $c(1-H_\infty(0))>1$, then $x_0$,
the smallest value in the support of $F_\infty$,
is positive, and $s_n^{(p)} \to x_0$ with probability $1$.
The value $x_0$ is the maximum of the function $z_\infty(m)$ for $m\in \mathbb{R}_+$.
\item [(ii)]
If $c(1-H_\infty(0))\le 1$ or $c(1-H_\infty(0))>1$ but $[a,b]$ is
not contained in $[0,x_0]$, then $m_\infty(b)<0$ and
\begin{equation}
  \mathbb{P} \bigl( \text{$s_{i_p}^{(p)}>b$ and $s_{i_p+1}^{(p)}<a$ for all large $p$}
\bigr) =1
\end{equation}
with $i_p$ defined in \eqref{eq:mainpropeqcon}.
(Here $s^{(p)}_0:= \infty$.)
\end{itemize}
\end{prop}

\section{Determination of $\supp(F_p)$}\label{sec:support}

The key part in applying Proposition \ref{prop:mainprop} turns out
to be determining the support of $F_p$. This can be extracted from
the following result due to Silverstein and Choi.

\begin{lem}[\cite{SilversteinC95}; see also Lemma 1.3 \cite{BaiS99}]
\label{lem:supportcriteria} If $x\notin \supp(F_p)$, then
$m:=m_p(x)$ satisfies
\begin{itemize}
\item[(i)] $m\in \mathbb{R}\setminus \{0\}$
\item[(ii)] $-\frac1{m} \notin \supp(H_p)$
\item[(iii)] $z'_p(m)>0$.
\end{itemize}
Conversely, if $m$ satisfies (i)-(iii), then $x=z_p(m)\notin
\supp(F_p)$.
\end{lem}

\begin{rem}
Lemma 1.3 of \cite{BaiS99} is stated for $H_\infty$. But the proof
of Lemma 1.3 in \cite{SilversteinC95} applies also to the finite
$p$ case of $H_p$ without any change. Indeed, the proposition
applies to any distribution defined by its Stieltjes transform as
in \eqref{eq:defFpviamp}.
\end{rem}

\begin{rem}
It is also shown in \cite{SilversteinC95} that $F_p$ has
continuous density on $\mathbb{R}_+$.
\end{rem}

When $T_p$ is as in \eqref{eq:Tpdef},
\begin{equation}\label{eq:Hp}
  dH_p(x)
= \frac1{p} \sum_{j=1}^{M} k_j \delta_{\alpha_j}(x) +
\biggl(1-\frac{r}{p}\biggr)\delta_1(x)
\end{equation}
and
\begin{equation}\label{eq:zp}
  z_p(m)
= -\frac1{m} + \frac{c_p}{1+m} + \frac1{n} \biggl[ \sum_{j=1}^M
\frac{k_j\alpha_j}{1+\alpha_j m} - \frac{r}{1+m} \biggr],
\end{equation}
where we recall that $r=k_1+\dots + k_M$. We first determine the
set of real $m$ such that $z_p'(m)>0$.

Now
\begin{equation}\label{eq:zpprime}
\begin{split}
  z_p'(m) &= \frac1{m^2} - \frac{c_p}{(1+m)^2} + \frac1{n} \biggl[
\sum_{j=1}^M \frac{-k_j\alpha_j^2}{(1+\alpha_j m)^2}
+\frac{r}{(1+m)^2} \biggr] \\
&= \frac{f(m) + \frac1n g(m)}{m^2(1+m)^2 \prod_{\ell=1}^M (1+\alpha_\ell m)^2},
\end{split}
\end{equation}
where
\begin{equation}\label{eq:fdef}
  f(m) := \bigl( (1+m)^2 -c_p m^2 \bigr) \prod_{\ell=1}^M (1+\alpha_\ell m)^2
\end{equation}
and
\begin{equation}\label{eq:gdef}
  g(m) := \biggl[
\sum_{j=1}^M \frac{-k_j\alpha_j^2}{(1+\alpha_j m)^2}
+\frac{r}{(1+m)^2} \biggr] m^2(1+m)^2  \prod_{\ell=1}^M (1+\alpha_\ell m)^2.
\end{equation}
We need the following basic lemmas of complex variables to determine
the solution of $z_p'(m)=0$.

\begin{lem}\label{lem:pert}
Let $h(z)$ be an analytic function in a closed disk
$\overline{D(z_0, r)}$ of radius $r>0$ centered at $z_0$. Then
there is $\epsilon_0>0$ such that for $0\le \epsilon<\epsilon_0$,
the equation
\begin{equation}\label{eq:perteq}
  z-z_0 = \epsilon h(z)
\end{equation}
has a unique solution in $D(z_0, r)$, which satisfies
\begin{equation}\label{eq:solzpert}
  z = z_0 + \epsilon h(z_0) + O(\epsilon^2).
\end{equation}
Furthermore, if $z_0$ is real and $h(z)$ is real for real $z$, the
solution \eqref{eq:solzpert} is real.
\end{lem}

\begin{proof}
As $h$ is continuous, there is a constant $C>0$ such that
$|h(z)|\le C$ for $|z-z_0|\le r$. When $|\epsilon| < \frac{r}{C}$,
for $|z-z_0|=r$,
\begin{equation}
  |z-z_0| =r>  |\epsilon| C \ge| \epsilon h(z)|.
\end{equation}
Hence from Rouche's theorem,
the number of zeros of $z-z_0-\epsilon h(z)$ inside $D(z_0, r)$
is equal to the number of zeros of $z-z_0$ inside $D(z_0,r)$, which is one.
The zero $z_\epsilon$ satisfies $z_\epsilon -z_0=\epsilon h(z_\epsilon) = O(\epsilon)$.
Thus
\begin{equation}
  z_\epsilon-z_0-\epsilon h(z_0) = \epsilon (h(z_\epsilon)-h(z_0))
= O(\epsilon^2).
\end{equation}
If $h(z)$ is real for real $z$, then by taking complex conjugate
of \eqref{eq:perteq}, we find that $\overline{z_\epsilon}$ is also
a solution. Since there is only one solution, we find that
$z_\epsilon$ is real.
\end{proof}

\begin{lem}\label{lem:pert2}
Let $h(z)$ be an analytic function in a closed disk
$\overline{D(z_0, r)}$ of radius $r>0$ centered at $z_0$ such that
$h(z_0)\neq 0$. Then there are $0<r_0\le r$ and $\epsilon_0>0$
such that for $0\le \epsilon<\epsilon_0$, the equation
\begin{equation}\label{eq:pert2eq}
  (z-z_0)^2 = \epsilon h(z)
\end{equation}
have precisely two distinct solutions in $D(z_0, r_0)$, which satisfy
\begin{equation}\label{eq:pert2z0}
  z = z_0 \pm \sqrt{\epsilon} \sqrt{h(z_0)} + O(\epsilon)
\end{equation}
where $\sqrt{h(z_0)}$ is an arbitrary branch.
Furthermore, suppose that
$z_0$ is real and $h(z)$ is real for real $z$.
Then if $h(z_0)>0$, both solutions \eqref{eq:pert2z0} are real.
On the other hand, if $h(z_0)<0$, both solutions \eqref{eq:pert2z0} are non-real.
\end{lem}

\begin{proof}
The proof of \eqref{eq:pert2z0} follows from Lemma \ref{lem:pert}
by taking the square root of \eqref{eq:pert2eq}.
When $z_0$ is real and $h(z)$ is real for real $z$,
the complex conjugate of a solution of \eqref{eq:pert2eq}
is also a solution.
Thus the two solutions\eqref{eq:pert2z0} of \eqref{eq:pert2eq}
are either complex conjugates of each other or both real since
there are precisely two distinct solutions.
Hence the Lemma follows.
\end{proof}

For the remainder of this section, we assume that $c\neq 1$ and none of
$\alpha_j$'s are equal to $1\pm \sqrt{c}$.
We further assume that $p$ and $n$ are sufficiently large so that
$c_p\neq 1$ and none of $\alpha_j$'s are equal to $1\pm \sqrt{c_p}$.
Then the numerator of \eqref{eq:zpprime} is a polynomial of degree exactly $2M+2$, and
we now determine all the solutions of $z_p'(m)=0$.

For $f$ defined in \eqref{eq:fdef}, the equation $f(m)=0$ has distinct
solutions
\begin{equation}\label{eq:mpm}
  m= \frac{-1}{1+\sqrt{c_p}}=:m_+, \qquad m=  \frac{-1}{1-\sqrt{c_p}}=:m_-
\end{equation}
of multiplicity $1$ and
\begin{equation}
  m= \frac{-1}{\alpha_j}, \qquad j=1,2,\dots, M,
\end{equation}
of multiplicity $2$. The roots of $z_p'(m)$ are expected to be
perturbations of the roots of $f(m)$, which we will find. First
consider $m_+$. Dividing the equation $f(m)+\frac1n g(m)=0$ by
$\frac{f(m)}{m-m_+}$, we obtain the equation
\begin{equation}
  m-m_+ + \frac1{n} \frac{m^2(1+m)^2}{(1-c_p)(m-m_-)} \biggl[
\sum_{j=1}^M \frac{-k_j\alpha_j^2}{(1+\alpha_j m)^2}
+\frac{r}{(1+m)^2} \biggr] =0.
\end{equation}
Lemma \ref{lem:pert} implies that there is a solution of
$z_p'(m)=0$ of the form
\begin{equation}
  m= m_+ + O\biggl( \frac1{n} \biggr),
\end{equation}
which is real. Similarly, there is a real solution of $z_p'(m)=0$
of the form
\begin{equation}
  m= m_- + O\biggl( \frac1{n} \biggr).
\end{equation}

Now consider the root $m=\frac{-1}{\alpha_j}$ of $f(m)=0$.
Dividing $f+\frac{1}{n} g=0$ by $\frac{f(m)}{(m+\frac1{\alpha_j})^2}$,
we obtain the equation
\begin{equation}\label{eq:malpha2}
  \biggl(m+\frac{1}{\alpha_j}\biggr)^2 = \frac1{n} G_j(m)
\end{equation}
where
\begin{equation}
  G_j(m) = \frac{-(1+\alpha_jm)^2m^2(1+m)^2}{\alpha_j^2(1-c_p)(m-m_+)(m-m_-)}
\biggl[
\sum_{\ell=1}^M \frac{-k_\ell\alpha_\ell^2}{(1+\alpha_\ell m)^2}
+\frac{r}{(1+m)^2} \biggr].
\end{equation}
Note that
\begin{equation}
  G_j\bigl(-\frac{1}{\alpha_j}\bigr)
=
\frac{k_j(\alpha_j-1)^2}{\alpha_j^4(1-c_p)(\frac{-1}{\alpha_j}-m_+)(\frac{-1}{\alpha_j}-m_-)}
\end{equation}
is not zero and also $G_j(m)$ is real for real $m$. Thus Lemma
\ref{lem:pert2} implies that there are precisely two solutions of
$z_p'(m)=0$  of the form
\begin{equation}\label{eq:msolalpha}
  m = -\frac{1}{\alpha_j} \pm \frac1{\sqrt{n}}
\sqrt{G_j\bigl(-\frac{1}{\alpha_j}\bigr)} + O\biggl( \frac1{n}
\biggr), \qquad j=1,\dots, M,
\end{equation}
where the pair for each $j$ are either both real or both non-real
depending on the sign of
$G_j\bigl(-\frac{1}{\alpha_j}\bigr)$.

Now when $c_p<1$, the condition
$G_j\bigl(-\frac{1}{\alpha_j}\bigr) >0$ is equivalent to
\begin{equation}
    \frac{-1}{\alpha_j} > m_+
\quad \text{or} \quad \frac{-1}{\alpha_j} < m_-,
\end{equation}
which is the same as
\begin{equation}
   \alpha_j >1+\sqrt{c_p} \quad \text{or}  \quad \alpha_j < 1-\sqrt{c_p}.
\end{equation}
On the other hand, when $c_p>1$, we note that $m_+<0<m_-$. The
condition $G_j\bigl(-\frac{1}{\alpha_j}\bigr) >0$ is now equivalent to
\begin{equation}
    m_+ < \frac{-1}{\alpha_j} < m_-,
\end{equation}
which is the same as (since $\alpha_j>0$)
\begin{equation}
   \alpha_j >1+\sqrt{c_p}.
\end{equation}

We summarize the above calculations.

\begin{lem}\label{lem:rootsofzp}
The solutions of $z_p'(m)=0$ are
\begin{equation}\label{eq:rootsofzp1}
    m= -\frac1{1+\sqrt{c_p}} + O\biggl( \frac1{n} \biggr)=:m_+^{(n)},
\qquad  m= - \frac1{1-\sqrt{c_p}} + O\biggl( \frac1{n} \biggr)=:m_-^{(n)}.
\end{equation}
and
\begin{equation}\label{eq:msolal}
  m = -\frac{1}{\alpha_j} \pm \frac1{\sqrt{n}}
\sqrt{G_j\bigl(-\frac{1}{\alpha_j}\bigr)}
+ O\biggl( \frac1{n} \biggr) =:m_{j, \pm}^{(n)},
\qquad j=1,\dots, M,
\end{equation}
all of multiplicity $1$.
Furthermore, the following holds.
\begin{itemize}
\item When $c_p<1$, $m_-^{(n)}<m_+^{(n)}<0$, and $m_{j,\pm}^{(n)}$
are real if and only if $\alpha_j>1+\sqrt{c_p}$ or
$\alpha_j<1-\sqrt{c_p}$. If $1-\sqrt{c_p}<\alpha_j <
1+\sqrt{c_p}$, $m_{j,+}^{(n)}$ and $m_{j, -}^{(n)}$ are complex
conjugates of each other.
\item When $c_p>1$,
$m_+^{(n)}<0<m_-^{(n)}$, and $m_{j,\pm}^{(n)}$ are real if and
only if $\alpha_j>1+\sqrt{c_p}$. If $\alpha_j<1+\sqrt{c_p}$,
$m_{j,+}^{(n)}$ and $m_{j, -}^{(n)}$ are complex conjugates of each
other.
\end{itemize}
\end{lem}

We now consider the cases when $c<1$ and when $c>1$ separately.

\subsection{When $c<1$}\label{sec:csmall1}


Let the indices $0\le M_0, M_1 \le M$ be defined as in Theorem
\ref{thm1} (recall that we assume that none of the $\alpha_j$'s are
equal to $1\pm \sqrt{c}$), so that
\begin{equation}\label{eq:M0def}
  \alpha_1>\dots >\alpha_{M_0} > 1+\sqrt{c} > \alpha_{M_0+1} >\dots >
\alpha_{M-M_1} > 1-\sqrt{c} > \alpha_{M-M_1+1} > \dots > \alpha_M.
\end{equation}
We now find the intervals in which $z_p'(m)>0$.

The denominator of \eqref{eq:zpprime} is non-negative.
>From Lemma \ref{lem:rootsofzp}, the numerator of \eqref{eq:zpprime} is factored as
\begin{equation}\label{eq:numerator}
   const \cdot (m-m_-^{(n)})(m-m_+^{(n)})
\prod_{j=1}^{M} (m-m_{j,-}^{(n)})  (m-m_{j,+}^{(n)})
\end{equation}
The constant prefactor is, from \eqref{eq:fdef} and \eqref{eq:gdef},
\begin{equation}\label{eq:prefactor}
  (1-c_p) \prod_{j=1}^M \alpha_j^2 + O\biggl( \frac1{n} \biggr),
\end{equation}
which is positive when $n$ is large enough. On the other hand,
among the terms in the product of \eqref{eq:numerator},
$m_{j,\pm}^{(n)}$ corresponding to the indices $M_0+1\le j\le M_1$
are complex conjugates of each other. Thus
\begin{equation}\label{eq:positiveproduct}
  \prod_{j=M_0+1}^{M_1+1} (m-m_{j,-}^{(n)})(m-m_{j,+}^{(n)}) \ge 0.
\end{equation}
Hence using the fact that
\begin{equation}
  0> m_{1,+}^{(n)} > m_{1,-}^{(n)} > \dots
> m_{M_0,+}^{(n)} > m_{M_0, -}^{(n)} > m_+^{(n)}
\end{equation}
and
\begin{equation}
  m_-^{(n)}> m_{M_1+1,+}^{(n)} > m_{M_1+1,-}^{(n)} > \dots
> m_{M,+}^{(n)} > m_{M, -}^{(n)},
\end{equation}
we find that the numerator of \eqref{eq:zpprime} is positve in the intervals
\begin{equation}\label{eq:pos1}
  (-\infty, m_{M,-}^{(n)})\cup (m_{M,+}^{(n)}, m_{M-1, -}^{(n)})\cup \dots \cup
 (m_{M_1+2,+}^{(n)}, m_{M_1+1, -}^{(n)})\cup (m_{M_1+1, +}^{(n)}, m_-^{(n)})
\end{equation}
union
\begin{equation}\label{eq:pos2}
   (m_+^{(n)}, m_{M_0,-}^{(n)}) \cup
(m_{M_0,+}^{(n)}, m_{M_0-1,-}^{(n)})
\cup \dots
\cup (m_{2,+}^{(n)},  m_{1,-}^{(n)}) \cup (m_{1,+}^{(n)}, \infty).
\end{equation}
The singular points of \eqref{eq:zpprime}
are not contained in any of the above intervals except for the
singular point $m=0$.
Hence the set of $m$ such that $z_p'(m)>0$ is equal to
\eqref{eq:pos1} union
\begin{equation}\label{eq:pos3}
   (m_+^{(n)}, m_{M_0,-}^{(n)}) \cup
(m_{M_0,+}^{(n)}, m_{M_0-1,-}^{(n)})
\cup \dots
\cup (m_{2,+}^{(n)},  m_{1,-}^{(n)}) \cup (m_{1,+}^{(n)}, 0)\cup (0,\infty).
\end{equation}

Now Lemma \ref{lem:supportcriteria} determines $\supp(F_p)$.

\begin{prop}\label{prop:suppFp}
Suppose that $c<1$ and none of $\alpha_j$ is equal to $1\pm
\sqrt{c}$. With the indices $M_0$ and $M_1$ defined in Theorem \ref{thm1},
for $n$ sufficiently large,
\begin{equation}\label{eq:suppFp}
\begin{split}
 \supp(F_p)^c
=  & (-\infty,0)\cup (0, z_{M,-}^{(n)}) \cup (z_{M,+}^{(n)}, z_{M-1, -}^{(n)})
\cup \dots \cup (z_{M_1+1, +}^{(n)}, z_-^{(n)}) \\
& \qquad \cup  (z_+^{(n)}, z_{M_0,-}^{(n)}) \cup
(z_{M_0,+}^{(n)}, z_{M_0-1,-}^{(n)})
\cup \dots
\cup (z_{2,+}^{(n)},  z_{1,-}^{(n)}) \cup (z_{1,+}^{(n)},\infty)
\end{split}
\end{equation}
where
\begin{equation}
   z_{\pm}^{(n)} = (1\pm \sqrt{c_p})^2 + O\biggl( \frac1{n} \biggr)
\end{equation}
and
\begin{equation}
   z_{j, \pm}^{(n)} = \alpha_j+\frac{c_p \alpha_j}{\alpha_j-1}
\pm \frac{A_j}{\sqrt{n}} + O\biggl( \frac{1}{n} \biggr), \qquad
j=1,\dots, M_0, \quad j=M_1+1, \dots, M,
\end{equation}
for some constant $A_j>0$.
The intervals in \eqref{eq:suppFp}
are disjoint.
\end{prop}

\begin{proof}
We will first see that the intervals \eqref{eq:pos1} union \eqref{eq:pos3}
satisfy the conditions (i)-(iii) of Lemma \ref{lem:supportcriteria}.
The condition (iii) is clearly satisfied.
Also $0$ is not contained in \eqref{eq:pos1} and \eqref{eq:pos3},
and so condition (i) is fullfilled.
Finally, as $\supp(H_p)=\{\alpha_1, \dots, \alpha_M, 1\}$
and
\begin{equation}
  m_-^{(n)}<-1<m_+^{(n)}, \qquad
m_{j,-}^{(n)}< -\frac1{\alpha_j} < m_{j,+}^{(n)},
\end{equation}
the condition (ii) is satisfied for $m$ in \eqref{eq:pos1} union \eqref{eq:pos3}.

We now need to find the image of the above intervals under
$z_p$.
Clearly, $z_p(-\infty)=0$, $z_p(0-)=+\infty$,
$z_p(0+)=-\infty$ and $z_p(+\infty)=0$.
A direct computation yields
\begin{equation}
  z_p( m_{\pm}^{(n)}) = (1\pm \sqrt{c_p})^2 + O\biggl( \frac1{n} \biggr).
\end{equation}
and
\begin{equation}
  z_p( m_{j, \pm}^{(n)})
= \alpha_j+\frac{c_p \alpha_j}{\alpha_j-1}
\pm \frac{A_j}{\sqrt{n}} + O\biggl( \frac{1}{n} \biggr)
\end{equation}
where
\begin{equation}
  A_j = \frac1{C_j} \biggl\{C_j^2 \alpha_j^2
\biggl( 1- \frac{c_p}{(\alpha_j-1)^2} \biggr)
+ k_j \biggr\}, \qquad C_j := \sqrt{G(-1/\alpha_j)}.
\end{equation}
Note that $A_j>0$ for $1\le j\le M_0$ and
$M_1+1\le j\le M$ since $\alpha_j>1+\sqrt{c_p}$ or
$\alpha_j<1-\sqrt{c_p}$.
Also it is straightforward to check from the graph of
the function
\begin{equation}
   x+ \frac{c_p x}{x-1}
\end{equation}
that
\begin{equation}
\begin{split}
  & 0< \alpha_M + \frac{c_p \alpha_M}{\alpha_M-1} <
\dots < \alpha_{M_1+1} + \frac{c_p \alpha_{M_1+1}}{\alpha_{M_1+1}-1} <
(1-\sqrt{c_p})^2 \\
&\qquad < (1+\sqrt{c_p})^2 <
\alpha_{M_0} + \frac{c_p \alpha_{M_0}}{\alpha_{M_0}-1}
< \dots < \alpha_1+ \frac{c_p \alpha_1}{\alpha_1-1}.
\end{split}
\end{equation}
This implies the Proposition.
\end{proof}

\subsection{When $c>1$}\label{sec:clarge1}

This case is similar to the previous case when $c<1$. We indicate
only the difference.

We again assume that $p$ and $n$ are large enough so that the set
of $j$'s satisfying $\alpha_j>1+\sqrt{c_p}$ is the same as the set of
$j$'s satisfying $\alpha_j>1+\sqrt{c}$. Let the index $0\le M_0\le
M$ be defined, as in Theorem \ref{thm1}. We further assume that
none of $\alpha_j$ is equal to $1+\sqrt{c}$ so that
\begin{equation}\label{eq:M0def2}
  \alpha_{M_0} > 1+\sqrt{c} > \alpha_{M_0+1}.
\end{equation}

The denominator of \eqref{eq:zpprime} is non-negative and as
before,  the numerator of \eqref{eq:zpprime} is equal to
\eqref{eq:numerator}. But this time, the constant prefactor
\eqref{eq:prefactor} is negative when $n$ is large enough. Also as
in \eqref{eq:positiveproduct},
\begin{equation}
  \prod_{j=M_0+1}^{M} (m-m_{j,-}^{(n)})(m-m_{j,+}^{(n)}) \ge 0.
\end{equation}
Now using the fact that
\begin{equation}
  m_-^{(n)}>0 > m_{1,+}^{(n)} > m_{1,-}^{(n)} > \dots
> m_{M_0,+}^{(n)} > m_{M_0, -}^{(n)} > m_+^{(n)},
\end{equation}
we find that the numerator of \eqref{eq:zpprime} is positive in
the intervals
\begin{equation}\label{eq:poscbig}
   (m_+^{(n)}, m_{M_0,-}^{(n)}) \cup
(m_{M_0,+}^{(n)}, m_{M_0-1,-}^{(n)}) \cup \dots \cup
(m_{2,+}^{(n)},  m_{1,-}^{(n)}) \cup (m_{1,+}^{(n)}, m_-^{(n)}).
\end{equation}
Hence taking into accounts of the singular point $m=0$ of
$z_p'(m)$, the intervals where $z_p'(m)>0$ is
\begin{equation}\label{eq:poscbig1}
   (m_+^{(n)}, m_{M_0,-}^{(n)}) \cup
(m_{M_0,+}^{(n)}, m_{M_0-1,-}^{(n)}) \cup \dots \cup
(m_{2,+}^{(n)},  m_{1,-}^{(n)}) \cup (m_{1,+}^{(n)}, 0) \cup (0,
m_-^{(n)}).
\end{equation}

The proof of the following proposition is parallel to Proposition
\ref{prop:suppFp}.

\begin{prop}\label{prop:suppFp2}
Suppose that $c>1$ and none of $\alpha_j$ is equal to $1 +
\sqrt{c}$. With the index $M_0$ defined in Theorem \ref{thm2}, for $n$ sufficiently large,
\begin{equation}\label{eq:suppFp2}
\begin{split}
 \supp(F_p)^c
=  & (-\infty,z_-^{(n)})\cup (z_+^{(n)}, z_{M_0,-}^{(n)}) \cup
(z_{M_0,+}^{(n)}, z_{M_0-1,-}^{(n)}) \cup \dots \cup
(z_{2,+}^{(n)},  z_{1,-}^{(n)}) \cup (z_{1,+}^{(n)},\infty)
\end{split}
\end{equation}
where
\begin{equation}
   z_{\pm}^{(n)} = (1\pm \sqrt{c_p})^2 + O\biggl( \frac1{n} \biggr)
\end{equation}
and
\begin{equation}
   z_{j, \pm}^{(n)} = \alpha_j+\frac{c_p \alpha_j}{\alpha_j-1}
\pm \frac{A_j}{\sqrt{n}} + O\biggl( \frac{1}{n} \biggr), \qquad
j=1,\dots, M_0,
\end{equation}
for some constant $A_j>0$.
The intervals in \eqref{eq:suppFp2} are disjoint.
\end{prop}

\section{Proof of Theorems \ref{thm1}, \ref{thm2} and
\ref{thm3}}\label{sec:proof}


When $T_p$ is \eqref{eq:Tpdef}, as $H_p$ is equal to
\eqref{eq:Hp},
\begin{equation}
  dH_\infty(x) = \delta_1(x).
\end{equation}
Hence
\begin{equation}\label{eq:zinf}
  z_{\infty}(m) = -\frac1{m} + c \int \frac{t}{1+tm} dH_\infty(t)
= -\frac1{m} + \frac{c}{1+m}.
\end{equation}
It is well-known that (\cite{SB95}, see also Theorem 3.4 of
\cite{BaiReview}) in this case,
\begin{equation}\label{eq:Finfformula}
  dF_\infty(\lambda) =
\begin{cases}
\frac1{2\pi \lambda}
\sqrt{((1+\sqrt{c})^2-\lambda)(\lambda-(1-\sqrt{c})^2)}
1_{[(1-\sqrt{c})^2, (1+\sqrt{c})^2]}(\lambda), \quad &c> 1 \\
\frac1{2\pi \lambda}
\sqrt{((1+\sqrt{c})^2-\lambda)(\lambda-(1-\sqrt{c})^2)}
1_{[(1-\sqrt{c})^2, (1+\sqrt{c})^2]}(\lambda) + (1-c) \delta_0 ,
\qquad &0<c\le 1.
\end{cases}
\end{equation}

\subsection{When $c<1$}

We first assume that none of $\alpha_j$ is equal to $1\pm
\sqrt{c}$ so that Proposition \ref{prop:suppFp} is applicable. The
case when some of $\alpha_j$ are equal to $1\pm \sqrt{c}$ will be
discussed at the end of this subsection.

When $T_p$ is \eqref{eq:Tpdef}, all the conditions (a)-(e) of
Proposition \ref{prop:mainprop} are satisfied or are defined
accordingly.

Now suppose $[a,b]$ is an interval satisfying condition (f). Since
\begin{equation}
  z_+^{(n)} \to (1+\sqrt{c})^2, \qquad z_-^{(n)}\to
  (1-\sqrt{c})^2,
\end{equation}
and for any $i$,
\begin{equation}\label{eq:zipm}
  z_{i,+}^{(n)}, \quad z_{i,-}^{(n)} \to \alpha_i + \frac{c \alpha_i}{\alpha_i-1},
\end{equation}
we see that
\begin{equation}
\begin{split}
  [a,b]\subset
  & \bigl(-\infty, 0 \bigr) \cup \biggl(0,  \alpha_M + \frac{c
  \alpha_M}{\alpha_M-1}\biggr) \cup
  ( \alpha_{M} + \frac{c \alpha_{M}}{\alpha_{M}-1},
  \alpha_{M-1} + \frac{c \alpha_{M-1}}{\alpha_{M-1}-1} \biggr) \\
  &\cup \dots
  \cup \biggl( \alpha_{M_1+1} + \frac{c \alpha_{M_1+1}}{\alpha_{M_1+1}-1},
  (1-\sqrt{c})^2 \biggr)  \\
  &\cup \biggl((1+\sqrt{c})^2,
\alpha_{M_0} + \frac{c \alpha_{M_0}}{\alpha_{M_0}-1} \biggr) \cup
\dots \cup \biggl( \alpha_2 + \frac{c \alpha_2}{\alpha_2-1},
\alpha_1 + \frac{c \alpha_1}{\alpha_1-1} \biggr) \cup \biggl(
\alpha_1 + \frac{c \alpha_1}{\alpha_1-1}, \infty\biggr).
\end{split}
\end{equation}
On the other hand,
\begin{equation}
  \supp(F_\infty)^c = (-\infty, 0) \cup \bigl(0, (1-\sqrt{c})^2\bigr)
  \cup \bigl( (1+\sqrt{c})^2, \infty \bigr).
\end{equation}
Hence $[a,b]\subset \supp(F_\infty)^c$. Also from definition
\eqref{eq:minf}, it is easy to see that $m'_\infty(z)>0$ for $z\in
\supp(F_\infty)^c$. The first consequence of (ii) of
Proposition \ref{prop:mainprop} (note that $H_\infty(0)=0$) is that
$m_\infty(b)<0$. Thus $m_\infty(a)< m_\infty(b) <0$.
Therefore, the condition \eqref{eq:mainpropeqcon} is equivalent
to the condition
\begin{equation}\label{eq:mainpropeqcondifferent}
 [a,b] \subset
 \bigl[z_\infty\bigl( -1/\sigma_{i_p+1}^{(p)} \bigr),
z_\infty\bigl( -1/\sigma_{i_p}^{(p)} \bigr) \bigr].
\end{equation}

We will consider four different choices of $[a,b]$.
First fix $1\le j\le M_0$.
Take
\begin{equation}
  [a,b] = \bigl[\alpha_j + \frac{c \alpha_j}{\alpha_j-1}+\epsilon,
\alpha_{j-1} + \frac{c \alpha_{j-1}}{\alpha_{j-1}-1}-\epsilon \bigr]
\end{equation}
for an arbitrary fixed $\epsilon >0$. (Here $\alpha_{0} :=
+\infty$.) From \eqref{eq:zipm}, we see that
\begin{equation}
  [a,b] \subset (z_{j, +}^{(n)}, z_{j-1,-}^{(n)})
\end{equation}
for all large $p$, and hence condition (f) is satisfied using
Proposition \ref{prop:suppFp}. Set
\begin{equation}
  i_p:= k_1+\dots + k_{j-1}.
\end{equation}
(When $j=1$, $i_p:=0$.)
For $T_p$ given by \eqref{eq:Tpdef},
\begin{equation}
  \sigma_{i_p}^{(p)}= \alpha_{j-1}, \qquad \sigma_{i_p+1}^{(p)}= \alpha_{j}.
\end{equation}
But
\begin{equation}
   z_\infty(-1/{\alpha_j}) = \alpha_j + \frac{c \alpha_j}{\alpha_j-1}
\end{equation}
and hence the condition \eqref{eq:mainpropeqcondifferent} is satisfied.
Therefore $i_p$ is defined to satisfy the condition \eqref{eq:mainpropeqcon}.
Proposition \ref{prop:suppFp} now implies that
\begin{equation}\label{eq:n413}
  \mathbb{P} \biggl(
\text{$s^{(p)}_{k_1+\dots +k_{j-1}} >
\alpha_{j-1} + \frac{c \alpha_{j-1}}{\alpha_{j-1}-1}-\epsilon$
and $s^{(p)}_{k_1+\dots + k_{j-1} +1}<
\alpha_j + \frac{c \alpha_j}{\alpha_j-1} + \epsilon$
for all large $p$} \biggr) =1.
\end{equation}
This yields that, $1\le j\le M_0-1$,
\begin{equation}
   \mathbb{P} \biggl(
\text{$\alpha_j+\frac{c \alpha_j}{\alpha_j-1} -\epsilon
< s_{k_1+\dots + k_{j-1}+k_j}^{(p)} \le \dots \le
s_{k_1+\dots k_{j-1}+1}^{(p)}
<  \alpha_j+\frac{c \alpha_j}{\alpha_j-1} +\epsilon$
for all large $p$} \biggr) =1,
\end{equation}
which implies \eqref{eq:thm11} for
$1\le j\le M_0-1$, and
\begin{equation}\label{eq:asymsM01}
  \mathbb{P} \biggl(
\text{$s_{k_1+\dots + k_{M_0-1}+1}^{(p)} <
\alpha_{M_0}+\frac{c \alpha_{M_0}}{\alpha_{M_0}-1} +\epsilon$
for all large $p$} \biggr) =1.
\end{equation}

For the second choice of $[a,b]$, set
\begin{equation}
  [a,b] = \bigl[(1+\sqrt{c})^2+\epsilon,
\alpha_{M_0} + \frac{c \alpha_{M_0}}{\alpha_{M_0}-1}-\epsilon \bigr]
\end{equation}
for an arbitrary fixed $\epsilon>0$.
Noting that
\begin{equation}
  z_+^{(n)} \to (1+\sqrt{c})^2
\end{equation}
and setting $i_p:= k_1+\dots + k_{M_0}$, a calculation similar to the above
yields that
\begin{equation}
  \mathbb{P} \biggl(
\text{$s^{(p)}_{k_1+\dots +k_{M_0}} >
\alpha_{M_0} + \frac{c \alpha_{M_0}}{\alpha_{M_0}-1}-\epsilon$
and $s^{(p)}_{k_1+\dots + k_{M_0} +1}<
(1+\sqrt{c})^2 + \epsilon$
for all large $p$} \biggr) =1.
\end{equation}
Thus, together with \eqref{eq:asymsM01}, we obtain
\eqref{eq:thm11} for $j=M_0$.
Also as \eqref{eq:zinf} and discussions around \eqref{eq:Flam}
implies that the support of the limiting spectral distribution of $B_p$
is $[(1-\sqrt{c})^2, (1+\sqrt{c})^2]$, we obtain \eqref{eq:thm12}.

As the third and forth choices of $[a,b]$, we set
\begin{equation}
  [a,b] = \bigl[\alpha_{M_1+1} + \frac{c \alpha_{M_1+1}}{\alpha_{M_1+1}-1} +\epsilon,
(1-\sqrt{c})^2-\epsilon \bigr]
\end{equation}
and
\begin{equation}
  [a,b] = \bigl[\alpha_{j+1} + \frac{c \alpha_{j+1}}{\alpha_{j+1}-1} +\epsilon,
\alpha_{j} + \frac{c \alpha_{j}}{\alpha_{j}-1}-\epsilon \bigr]
\end{equation}
for some $M_1+1\le j\le M$ ($\alpha_{M+1}:=0$), respectively.
Arguments as above imply the remaining part of Theorem \ref{thm1}.

\bigskip

We now consider the case when an $\alpha_j$ is equal to $1\pm
\sqrt{c}$. We first observe certain monotonicity of the
eigenvalues $s^{(p)}_j$ on $\alpha_j$s. Note that the matrix
$\underline{B}_p:= \frac1{n} Z_p' T_p Z_p$ has the same set of
eigenvalues as $B_p$ except for $|p-n|$ zero eigenvalues. Consider
a set of parameters $\beta_j$, $1\le j\le M$, such that
$\alpha_j\ge \beta_j$. Let $\hat{T}_p$ be the matrix $T_p$ with
$\alpha_j$'s replaced by $\beta_j$'s, and set $\hat{B}_p=
\frac1{n} \hat{T}_p^{1/2}Z_p Z_p' \hat{T}_p^{1/2}$ and
$\underline{\hat{B}}_p= \frac1{n} Z_p' \hat{T}_p Z_p$. Then
clearly, $\underline{B}_p$ and $\underline{\hat{B}}_p$ are
Hermitian, and $\underline{B}_p \ge \underline{\hat{B}}_p$. Hence
from the min-max principle (see e.g. \cite{HornJ}), we find that
\begin{equation}
  s^{(p)}_j \ge \hat{s}^{(p)}_j
\end{equation}
for all non-zero eigenvalues, where $\hat{s}^{(p)}_j$ denotes the
eigenvalues of $\hat{B}_p$.

Suppose that
\begin{equation}\label{eq:alphas}
  \alpha_1> \dots > \alpha_{M_0} > 1+\sqrt{c} = \alpha_{M_0+1}
> \dots > \alpha_{M-M_1} > 1-\sqrt{c} > \alpha_{M-M_1+1}
>\dots > \alpha_M.
\end{equation}
Replacing in \eqref{eq:alphas} $\alpha_{M_0+1}$ by
$(1\pm\epsilon)\alpha_{M_0+1}=(1\pm\epsilon)(1+\sqrt{c})$ for
sufficiently small $\epsilon>0$, the above monotonicity argument
implies the following:
\begin{itemize}
\item[(i)] For each $1\le j\le M_0$,
\begin{equation}
  \alpha_{j} + \frac{c
  \alpha_{j}}{\alpha_j-1}
  \le \liminf s^{(p)}_{k_1+\dots + k_{j-1}+i}
  \le \limsup s^{(p)}_{k_1+\dots + k_{j-1}+i}
  \le \alpha_{j} + \frac{c
  \alpha_{j}}{\alpha_j-1},
  \qquad 1\le i\le k_j.
\end{equation}
almost surely.
\item[(ii)]
\begin{equation}
  (1+\sqrt{c})^2 \le \liminf s^{(p)}_{k_1+\dots+k_{M_0}+1}
  \le \limsup s^{(p)}_{k_1+\dots+k_{M_0}+1}
  \le (1+\epsilon)\alpha_{M_0+1} + \frac{c (1+\epsilon)\alpha_{M_0+1}}{(1+\epsilon)\alpha_{M_0+1}-1}
\end{equation}
almost surely.
\item[(iii)]
\begin{equation}
  (1-\sqrt{c})^2 \le
  \liminf s^{(p)}_{p-r+k_1+\dots + k_{M_1}} \le
  \limsup s^{(p)}_{p-r+k_1+\dots + k_{M_1}} \le
  (1-\sqrt{c})^2
\end{equation}
almost surely.
\item[(iv)]
For each $M_1+1\le j\le M$,
\begin{equation}
\begin{split}
  & \alpha_{j} + \frac{c
  \alpha_{j}}{\alpha_j-1} \le
  \liminf s^{(p)}_{p-r+k_1+\dots + k_{j-1} +i}\to  \alpha_{j} + \frac{c
  \alpha_{j}}{\alpha_j-1} \\
  &\qquad
  \le  \limsup s^{(p)}_{p-r+k_1+\dots + k_{j-1} +i}\to  \alpha_{j} + \frac{c
  \alpha_{j}}{\alpha_j-1}
  \le \alpha_{j} + \frac{c
  \alpha_{j}}{\alpha_j-1}
  \qquad 1\le  i\le k_{j}
\end{split}
\end{equation}
almost surely.
\end{itemize}
Since
\begin{equation}
  \lim_{\epsilon\downarrow 0}
  (1+\epsilon)\alpha_{M_0+1} + \frac{c
  (1+\epsilon)\alpha_{M_0+1}}{(1+\epsilon)\alpha_{M_0+1}-1}
  = (1+\sqrt{c})^2
\end{equation}
and the above result is true for arbitrary sufficiently small
$\epsilon>0$, Theorem \ref{thm1} follows for the case when the
parameters are given by \eqref{eq:alphas}. For the case when
$\alpha_{M-M_1}=1-\sqrt{c}$, the argument is almost identical, and
we skip the details.

\subsection{When $c>1$}

>From \eqref{eq:Finfformula}, when $c>1$, the smallest value in the
support of $F_\infty$ is
\begin{equation}
 x_0=(1-\sqrt{c})^2 >0.
\end{equation}
Hence Proposition \ref{prop:mainprop} (i) implies that
\begin{equation}
  s_n^{(p)} \to (1-\sqrt{c})^2.
\end{equation}
Since when $p>n$, at least $p-n$ eigevalues $s^{(p)}_j$ are equal to $0$,
we conclude that
\begin{equation}
  s_{n+1}^{(p)} =\dots = s_p^{(p)} =0.
\end{equation}
Therefore, \eqref{eq:thm23} and \eqref{eq:thm24} are obtained.

The proof of \eqref{eq:thm21} and \eqref{eq:thm22} is similar to
the case when $c<1$ by using Proposition \ref{prop:suppFp2} and
noting that an interval $[a,b]$ satisfying condition (f) of
Proposition \ref{prop:mainprop} is contained in
\begin{equation}
\begin{split}
  \bigl(-\infty, (1-\sqrt{c})^2 \bigr)
  \cup \biggl((1+\sqrt{c})^2,  \alpha_{M_0} + \frac{c
  \alpha_{M_0}}{\alpha_{M_0}-1} \biggr) \cup \dots \cup
  \biggl( \alpha_2 + \frac{c
  \alpha_2}{\alpha_2-1}, \alpha_1 + \frac{c
  \alpha_1}{\alpha_1-1}\biggr)
  \cup \biggl( \alpha_1 + \frac{c
  \alpha_1}{\alpha_1-1}, \infty\biggr),
\end{split}
\end{equation}
which is a subset of
\begin{equation}
  \supp(F_\infty)^c = (-\infty, (1-\sqrt{c})^2\bigr)
  \cup \bigl( (1+\sqrt{c})^2, \infty \bigr).
\end{equation}

\subsection{When $c=1$}\label{sec:cis1}

Since the limiting distribution \eqref{eq:PM} for $c=1$ has a
continuous density on the interval $(0,4)$, it is easy to see
\eqref{eq:thm33}.

We first observe a monotonicity of $s^{(p)}_j$ in $n$. Let
$\hat{Z}_p = (Z_{ij})$, $1\le i\le p, 1\le j\le \hat{n}$ and let
$\hat{B}_p:= \frac1{\hat{n}} T_p^{1/2} \hat{Z}_p\hat{Z}_p'
T_p^{1/2}$. When $\hat{n}> n$, it is clear that
\begin{equation}
  \hat{n} \hat{B}_p \ge n B_p.
\end{equation}
Therefore, if the ordered eigenvalues of $\hat{B}_p$ are denoted
by $\hat{s}^{(p)}_p$, the min-max principle implies that
\begin{equation}\label{eq:monoton}
  \hat{n} \hat{s}^{(p)}_j \ge n s^{(p)}_j
\end{equation}
for all $1\le j\le p$.

Take $\hat{n}= \bigl[ \frac{n}{1+ \epsilon} \bigr]$ for
$\epsilon>0$ where $[x]$ denotes the largest integer $\le x$. Then
for sufficiently small $\epsilon>0$,
\begin{equation}
  \alpha_1>\dots > \alpha_{M_0} > 1+\sqrt{1+\epsilon}
  > \alpha_{M_0+1}>\dots >\alpha_M.
\end{equation}
By applying Theorem \ref{thm2} and using \eqref{eq:monoton}, we
obtain the following:
\begin{itemize}
\item
For each $1\le j\le M_0$,
\begin{equation}\label{eq:c11}
  \liminf s^{(p)}_{k_1+\dots + k_{j-1}+i} \ge \frac1{1+\epsilon} \biggl(
  \alpha_{j} + \frac{(1+\epsilon)\alpha_{j}}{\alpha_j-1} \biggr),
  \qquad 1\le i\le k_j.
\end{equation}
almost surely.
\item
\begin{equation}\label{eq:c12}
  \liminf s^{(p)}_{k_1+\dots +k_{M_0}+1} \ge \frac1{1+\epsilon}
  (1+\sqrt{1+\epsilon})^2
\end{equation}
almost surely.
\end{itemize}

On the other hand, take $\hat{n}= \bigl[ \frac{n}{1- \epsilon}
\bigr]$ for $\epsilon>0$. We first assume $2> \alpha_{M_0+1}$.
Then as $\alpha_M>0$, for sufficiently small $\epsilon>0$,
\begin{equation}
  \alpha_1>\dots > \alpha_{M_0} > 1+\sqrt{1-\epsilon}
  > \alpha_{M_0+1}>\dots >\alpha_M > 1-\sqrt{1-\epsilon},
\end{equation}
and hence $M_1=M$. By applying Theorem \ref{thm1} and using
\eqref{eq:monoton}, we obtain the following:
\begin{itemize}
\item
For each $1\le j\le M_0$,
\begin{equation}\label{eq:c111}
  \limsup s^{(p)}_{k_1+\dots + k_{j-1}+i} \le
  \frac1{1-\epsilon} \biggl(
  \alpha_{j} + \frac{(1-\epsilon)\alpha_{j}}{\alpha_j-1} \biggr),
  \qquad 1\le i\le k_j.
\end{equation}
almost surely.
\item
\begin{equation}\label{eq:c112}
  \limsup s^{(p)}_{k_1+\dots +k_{M_0}+1} \le
  \frac1{1-\epsilon}
  (1+\sqrt{1-\epsilon})^2
\end{equation}
almost surely.
\end{itemize}
If $\alpha_{M_0+1}=2$, then for sufficiently small $\epsilon>0$,
\begin{equation}
  \alpha_1>\dots > \alpha_{M_0} > \alpha_{M_0+1} > 1+\sqrt{1-\epsilon}
  > \alpha_{M_0+2}>\dots >\alpha_M > 1-\sqrt{1-\epsilon}.
\end{equation}
Hence Theorem \ref{thm1} implies \eqref{eq:c111} but
\eqref{eq:c112} becomes
\begin{equation}\label{eq:c113}
  \limsup s^{(p)}_{k_1+\dots +k_{M_0}+1} \le
  \frac1{1-\epsilon} \biggl(
  \alpha_{M_0+1} + \frac{(1-\epsilon)\alpha_{M_0+1}}{\alpha_{M_0+1}-1} \biggr)
\end{equation}
almost surely.

Therefore \eqref{eq:c11} and \eqref{eq:c111} yield
\eqref{eq:thm31}, and \eqref{eq:c12}, \eqref{eq:c112} and
\eqref{eq:c113} yield \eqref{eq:thm32}.



\end{document}